\documentclass[12pt,reqno]{amsart}
\usepackage{amsmath,amsthm,amssymb,amsfonts,amscd}
\usepackage{mathrsfs}
\usepackage{bbm}
\usepackage{booktabs}
\usepackage{float}
\usepackage[toc,page]{appendix}
\usepackage{makecell}
\usepackage{bbding}
\usepackage[backref=page]{hyperref}
\usepackage{geometry}
\usepackage{color}
\usepackage{xcolor}
\usepackage{picture,epic}
\usepackage{tikz}

\usepackage{etoolbox}
\usepackage{orcidlink}

\numberwithin{equation}{section}

\theoremstyle{plain}
\newtheorem{theorem}{Theorem}[section]

\newtheorem{proposition}[theorem]{Proposition}

\theoremstyle{definition}

\theoremstyle{remark}

\renewcommand{\Re}{\operatorname{Re}}

\newcommand{\propertyS}[1]{%
  \par\noindent\textbf{Property (S).} #1 \par
}

\makeatletter
\def\@tocline#1#2#3#4#5#6#7{\relax
  \ifnum #1>\c@tocdepth 
  \else
    \par \addpenalty\@secpenalty\addvspace{#2}%
    \begingroup \hyphenpenalty\@M
    \@ifempty{#4}{%
      \@tempdima\csname r@tocindent\number#1\endcsname\relax
    }{%
      \@tempdima#4\relax
    }%
    \parindent\z@ \leftskip#3\relax \advance\leftskip\@tempdima\relax
    \rightskip\@pnumwidth plus4em \parfillskip-\@pnumwidth
    #5\leavevmode\hskip-\@tempdima
      \ifcase #1
       \or\or \hskip 1em \or \hskip 2em \else \hskip 3em \fi%
      #6\nobreak\relax
    \hfill\hbox to\@pnumwidth{\@tocpagenum{#7}}\par
    \nobreak
    \endgroup
  \fi}
\makeatother

\allowdisplaybreaks

\begin{document}

\title{Discrete restrictions from Laurent monomial systems for multiple Dirichlet series}

\author{Shenghao Hua~\orcidlink{0000-0002-7210-2650}}

\subjclass[2020]{05A15,11M06,11A07}

\keywords{Multiple Dirichlet series}

\begin{abstract}
We introduce a special class of multiple Dirichlet series whose terms are supported on a variety and which admit an Euler product structure.
\end{abstract}

\maketitle
\section{Introduction}

The restriction of the Fourier or Mellin transform to submanifolds is an active area of research.
A Dirichlet series can be viewed as a discrete sampling and weighted summation of the Mellin transform at integer points, effectively replacing the continuous measure by a Stieltjes integral. There is a profound connection between the two: for example, finite partial sums of a Dirichlet series can be recovered via a special form of the inverse Mellin transform known as Perron's formula. Moreover, the asymptotic functional equation expansion of an automorphic \( L \)-function is obtained through the inverse Mellin transform of the completed \( L \)-function.
If we view a multiple Dirichlet series as a discrete Mellin transform integral over the entire Euclidean space, a natural question arises: \emph{can such an integral be restricted to integrals over submanifolds, such as algebraic varieties?}

For an algebraic variety \( V \) in the \( t \)-dimensional Euclidean space defined by \( m \leq t \) defined by a system of rational functional equations with integer coefficients
\[
f_1(\boldsymbol{x}) = 0, \dots, f_m(\boldsymbol{x}) = 0,
\]
where $\boldsymbol{x}=(x_1,\dots,x_t)$, we consider restricting the summation to its positive integer solutions.
For the multiple Dirichlet series
\[
\sum_{\substack{n_1,\dots,n_t \geq 1}}
\frac{ a(\boldsymbol{n})}
{\prod_{i=1}^t n_i^{s_i}},
\]
defined by coefficients \( a(\boldsymbol{n}) \), where $\boldsymbol{n}=(n_1,\dots,n_t)$,
the series converges when all \( \Re s_i \) are sufficiently large.
Here we consider a multiple Dirichlet series twisted by the indicator function of \( V \), which is of the form
\[
\sum_{\substack{n_1,\dots,n_t \geq 1}}
\frac{
\delta_{\boldsymbol{n} \in V(\mathbb{Z})} a(\boldsymbol{n})}{\prod_{i=1}^{t} n_i^{s_i}}
=
\sum_{\substack{n_1,\dots,n_t \geq 1 \\ \boldsymbol{n} \in V(\mathbb{Z})}}
\frac{a(\boldsymbol{n})}{\prod_{i=1}^{t} n_i^{s_i}},
\]
where \( \delta_{\boldsymbol{n} \in V(\mathbb{Z})} \) denotes the indicator function of the set of integer points on \( V \).

Such a summation is easily shown to converge when the real parts of the variables are sufficiently large. However, we hope that it also possesses desirable properties, similar to those enjoyed by functions in the Selberg class.
The Selberg class of functions was introduced by Selberg in the 1980s as an axiomatic framework to generalize and study properties of \( L \)-functions, such as the Riemann zeta function and Dirichlet \( L \)-functions. These functions satisfy a set of analytic conditions, including a functional equation, Euler product that is locally generated by finitely many terms, and Ramanujan hypothesis. Around the same time, the theory of automorphic \( L \)-functions developed rapidly as part of the Langlands program, which seeks deep connections between number theory and representation theory. Automorphic \( L \)-functions arise from automorphic forms and play a crucial role in modern number theory, encompassing many classical \( L \)-functions within a unified theory. The interplay between the Selberg class and automorphic \( L \)-functions continues to be a central theme in analytic number theory.
Kaczorowski and Perelli~\cite{KaczorowskiPerelli3}, under certain conditions, proved a conjecture of Sarnak concerning functions in the Selberg class. Specifically, they showed that it is possible to find a countable set of primitive Selberg class functions such that all Selberg class functions can be generated as products of shifts of these primitive elements.

We first present several basic properties of these multiple Dirichlet series and then offer a special (and possibly general) example. In the last section, we give a counterexample to a natural geometric characterization suggested by Property~(S).

\section{Multiple Dirichlet Series Restricted to a Variety}

We aim to investigate the relationship between Selberg class functions and automorphy, with particular emphasis on the significance of being locally generated by finitely many terms.
To this end, we first construct an Euler product by taking \( a(\boldsymbol{n}) \) to be a multiplicative function. However, this Euler product cannot be expressed locally by finitely many terms, which distinguishes it from functions in the Selberg class.
And the set of positive integer solutions to \( V \) must be closed under multiplication of prime power factors arising from these solutions.

We define the following property:

\propertyS{
Let
\[
(n_1 = \prod_p p^{\alpha_{1,p}}, \dots, n_t = \prod_p p^{\alpha_{t,p}}), \quad
(n_1' = \prod_p p^{\beta_{1,p}}, \dots, n_t' = \prod_p p^{\beta_{t,p}}) \in V(\mathbb{Z})
\]
be two integer points on the variety \( V \). Then, for every prime \( p \) and for every choice of exponents
\[
(\gamma_{1,p}, \dots, \gamma_{t,p}) \in \{ (\alpha_{1,p}, \dots, \alpha_{t,p}),\; (\beta_{1,p}, \dots, \beta_{t,p}) \},
\]
we have
\[
\left( \prod_p p^{\gamma_{1,p}}, \dots, \prod_p p^{\gamma_{t,p}} \right) \in V(\mathbb{Z}).
\]
}

Then we have
\begin{theorem}
Let $V$ be a variety. The twisted multiple Dirichlet series
\begin{equation}\label{eqn:wtDs}
\sum_{\substack{n_1,\dots,n_t \geq 1 \\ \boldsymbol{n} \in V(\mathbb{Z})}}
\frac{a(\boldsymbol{n})}{\prod_{i=1}^{t} n_i^{s_i}}
\end{equation}
can be represented as an Euler product for any given arithmetic function $a(\boldsymbol{\cdot})$ if and only if $V$ satisfies Property~(S).
\end{theorem}

\begin{proof}
The condition that $V$ satisfies Property~(S) is equivalent to the statement that the set of positive integer solutions to $V$ decomposes as a Cartesian product over primes:
\[
V(\mathbb{Z}_{>0}) = \prod_p V_p,
\]
where each \( V_p \subset \mathbb{Z}_{\geq 0}^t \) denotes the set of admissible exponent tuples at the prime $p$.

Suppose the series~\eqref{eqn:wtDs} admits an Euler product decomposition. Then for each prime $p$, we can identify a corresponding subset \( V_p \subset \mathbb{Z}_{\geq 0}^t \) capturing the admissible exponent tuples at $p$.

Conversely, if $V$ satisfies Property~(S), then using the sets \( V_p \) for each prime $p$, one can factor the sum in~\eqref{eqn:wtDs} into an Euler product accordingly.
\end{proof}

\section{An example}

We consider the following equivalent formulation of a system of Laurent monomial equations:
\[
f_i(x_1,\ldots,x_t) = 0 \quad \text{for} \quad i = 1,\ldots,m,
\]
where each \( f_i \) is of the form
\[
\omega_i f_i(x_1,\ldots,x_t) = \prod_{j=1}^{t} x_j^{a_{i,j}} - \omega'_i,
\]
with \( a_{i,j} \in \mathbb{Z} \), and \( \omega_i,\omega'_i \in \mathbb{Z}_{\geq 1} \), ensuring that the system admits potential positive solutions.

Let $A=(a_{i,j})_{\substack{1\leq i\leq m\\1\leq j\leq t}}$, $\boldsymbol{s}=(s_1,\dots,s_t)$,
 $\boldsymbol{\omega}=(\omega_1,\dots,\omega_m)$,
 $\boldsymbol{\omega'}=(\omega'_1,\dots,\omega'_m)$, and
\[
D_{A,a(\boldsymbol{\cdot})}
(\boldsymbol{s};\boldsymbol{\omega},\boldsymbol{\omega'})
:=
\sum_{\substack{n_1,\dots,n_t \geq 1 \\ \omega_i\prod_{j=1}^{t}n_j^{a_{i,j}}=\omega'_i,~1\leq i\leq m}}
\frac{a(\boldsymbol{n})}{\prod_{i=1}^{t} n_i^{s_i}}.
\]

When all $\Re(s_i)$ are sufficiently large, we have the expected Euler product representation:
\[
D_{A,a(\boldsymbol{\cdot})}
(\boldsymbol{s};\boldsymbol{\omega},\boldsymbol{\omega'}) =
\prod_{\substack{p~\text{prime}}}
\left(
\sum_{\substack{\alpha_1,\dots,\alpha_t \geq 0 \\
v_p(\omega'_i)-v_p(\omega_i)
+\sum_{j=1}^{t}a_{i,j}\alpha_j
=0 \\
\text{for } 1 \leq i \leq m}}
\frac{a(\boldsymbol{p^\alpha})}
{p^{\sum_{i=1}^{t}s_i\alpha_i}}
\right),
\]
where $v_p(x)$ is the $p$-adic valuation of $x$, and $\boldsymbol{p^\alpha}=(\alpha_1,\dots,\alpha_t)$.

According to the definition, we obtain
\[
D_{A, a(\boldsymbol{\cdot})}(\boldsymbol{s}; \boldsymbol{\omega},\boldsymbol{\omega'}) =
D_{-A, a(\boldsymbol{\cdot})}(\boldsymbol{s};
\boldsymbol{\omega'},\boldsymbol{\omega}).
\]

Let
\[
A^{(3)}= \begin{pmatrix} A^{(1)} & 0 \\ 0 & A^{(2)} \end{pmatrix},
\]
\[
\boldsymbol{s^{(3)}} = (\boldsymbol{s^{(1)}}, \boldsymbol{s^{(2)}}), \quad
\boldsymbol{\omega^{(3)}} = (\boldsymbol{\omega^{(1)}}, \boldsymbol{\omega^{(2)}}), \quad
\boldsymbol{\omega'^{(3)}} = (\boldsymbol{\omega'^{(1)}}, \boldsymbol{\omega'^{(2)}}),
\]
\[
\boldsymbol{n^{(3)}} = (\boldsymbol{n^{(1)}}, \boldsymbol{n^{(2)}}), \quad
a^{(3)}(\boldsymbol{n^{(3)}}) = a^{(1)}(\boldsymbol{n^{(1)}}) \cdot a^{(2)}(\boldsymbol{n^{(2)}}).
\]

By definition, we then have
\[
D_{A^{(1)}, a^{(1)}(\cdot)}(\boldsymbol{s^{(1)}}; \boldsymbol{\omega^{(1)}}, \boldsymbol{\omega'^{(1)}}) \cdot D_{A^{(2)}, a^{(2)}(\cdot)}(\boldsymbol{s^{(2)}}; \boldsymbol{\omega^{(2)}}, \boldsymbol{\omega'^{(2)}}) = D_{A^{(3)}, a^{(3)}(\cdot)}(\boldsymbol{s^{(3)}}; \boldsymbol{\omega^{(3)}}, \boldsymbol{\omega'^{(3)}}).
\]

Moreover, by definition, the following elementary row operations in Gaussian elimination do not affect the value of \( D_{A, a(\cdot)}(\boldsymbol{s}; \boldsymbol{\omega}, \boldsymbol{\omega'}) \):
\begin{itemize}
\item Swapping the \( i \)-th row of \( A \) with the \( j \)-th row.
  \item Multiplying the \( i \)-th row of \( A \) by \(-1\), and exchanging \(\omega_i\) and \(\omega'_i\).
  \item Adding \( b \) times the \( j \)-th row of \( A \) to the \( i \)-th row, where \( b \in \mathbb{Z} \), and replacing \(\omega_i\) by \(\omega_i \omega_j^b\), similarly for \(\omega'_i\).
\end{itemize}

In particular, let
\[
a(\boldsymbol{n}) = \prod_{j=1}^{t} \lambda_{\Pi_j}(n_j),
\]
where each \(\lambda_{\Pi_j}(\cdot)\) denotes the normalized Fourier coefficient of the \(\mathrm{GL}(k_j,\mathbb{R})\) automorphic form \(\Pi_j\).
In particular, when $\Pi_j = \mathrm{1}$, $L(s, \Pi_j)=\zeta(s)$.
Let $\boldsymbol{\Pi}=(\Pi_1,\dots,\Pi_m)$, we define
\[
D_{A, \boldsymbol{\Pi}}(\boldsymbol{s}; \boldsymbol{\omega}, \boldsymbol{\omega'}) := D_{A,a(\boldsymbol{\cdot})}(\boldsymbol{s}; \boldsymbol{\omega}, \boldsymbol{\omega'}).
\]

There is a natural origin: when all $\Re s_j > 1$, this series can be derived
from certain special twisted moments of a family of automorphic $L$-functions
associated with Dirichlet twists. Let $q$ be a sufficiently large prime.
Then there exists a constant $\eta > 0$, depending on the coefficients $a_{i,j}$
and the real parts $\Re s_j$, such that
\begin{multline*}
D_{A, \boldsymbol{\Pi}}(\boldsymbol{s}; \boldsymbol{\omega}, \boldsymbol{\omega'}) =
\frac{1}{(q-1)^{m}} \sum_{\substack{\chi_1, \dots, \chi_m \bmod{q}}}
\prod_{j=1}^{t} L\!\left(s_j, \Pi_j \times \prod_{i=1}^{m} \chi_i^{a_{i,j}} \right)
\prod_{i=1}^{m} \chi_i(\omega_i)\overline{\chi_i}(\omega'_i)
\\
+ O_{A,\boldsymbol{\Pi},\boldsymbol{s},\boldsymbol{\omega},\boldsymbol{\omega'}}
\!\left(q^{-\eta}\right).
\end{multline*}

\section{A Counterexample}

Property~(S) constrains only the positive integral points of a variety and, without an additional density hypothesis, does not determine its algebraic structure. In particular, the assertion that every variety satisfying Property~(S) is equivalent to a Laurent monomial system is false, even for an irreducible variety with infinitely many positive integral points.

\begin{proposition}\label{prop:conter}
There exists an irreducible variety \(V\) with infinitely many positive integral points such that \(V\) satisfies Property~(S), but \(V(\mathbb{Z}_{>0})\) is not the positive integral solution set of any Laurent monomial system.
\end{proposition}

\begin{proof}
Let \(V\subset \mathbb{A}^{3}_{\mathbb{Q}}\), with coordinates \((x,y,z)\), be defined by
\begin{equation}\label{eqn:counterexample}
10(x+y-5)^2+(x-y)^2=9,
\end{equation}
where \(z\) is free. After the affine change of variables
\[
u=x+y-5,\qquad v=x-y,
\]
the defining polynomial becomes \(10u^2+v^2-9\). Its homogenization
\[
10u^2+v^2-9w^2
\]
is a nondegenerate ternary quadratic form, so the corresponding plane conic is geometrically irreducible. It follows that \(V\) is irreducible.

If \(x,y\) are positive integers satisfying~\eqref{eqn:counterexample} and \(x+y-5\neq 0\), then the first term on the left-hand side is at least \(10\), which is impossible. Hence \(x+y=5\), and then \((x-y)^2=9\). Therefore
\begin{equation}\label{eqn:counterexample-points}
V(\mathbb{Z}_{>0})
=\{(1,4,n):n\geq 1\}\mathbin{\cup}\{(4,1,n):n\geq 1\}.
\end{equation}

We verify Property~(S). For every prime \(p\neq 2\), the \(p\)-adic valuation vector of the first two coordinates of either family in~\eqref{eqn:counterexample-points} is \((0,0)\). At \(p=2\), it is either \((0,2)\) or \((2,0)\). Thus any prime-by-prime recombination of two points in~\eqref{eqn:counterexample-points} again has first two coordinates equal to either \((1,4)\) or \((4,1)\); its third coordinate is a positive integer. The recombined point therefore still lies in \(V(\mathbb{Z}_{>0})\).

It remains to show that this positive integral solution set cannot arise from a Laurent monomial system. Suppose, to the contrary, that such a system exists. Write one of its equations as
\[
\omega_i x^{a_i}y^{b_i}z^{c_i}=\omega'_i,
\]
where \(a_i,b_i,c_i\in\mathbb{Z}\) and \(\omega_i,\omega'_i\in\mathbb{Z}_{\geq 1}\). Since the equation holds at \((1,4,n)\) for every \(n\geq 1\), comparison of \(n=1\) and \(n=2\) gives \(c_i=0\). Evaluating it at \((1,4,1)\) and \((4,1,1)\) then gives
\[
4^{b_i}=4^{a_i},
\]
and hence \(a_i=b_i\). Consequently, for every \(n\geq 1\),
\[
\omega_i 2^{a_i}2^{b_i}
=\omega_i4^{a_i}
=\omega'_i,
\]
so \((2,2,n)\) satisfies every equation in the system. On the other hand, substituting \((2,2,n)\) into the left-hand side of~\eqref{eqn:counterexample} gives \(10\), not \(9\), a contradiction.
\end{proof}

The obstruction is that Property~(S) only sees the set of positive integral points. A possible corrected formulation would at least require \(V(\mathbb{Z}_{>0})\) to be Zariski dense in \(V\), together with a precise meaning of equivalence, for example that \(V\cap\mathbb{G}_m^t\) is a translate of a subtorus or that its Laurent ideal is generated by binomials. We do not pursue this geometric question here, and turn instead to the two analytic conjectures.


\end{document}